\documentclass{amsart}

\usepackage{graphicx} 
\usepackage{amsmath,amssymb,amsthm,tikz}
\usepackage{xcolor}
\usepackage{esint}
\usepackage{pgfplots}
\usepackage{scalerel}
\pgfplotsset{compat=1.18}
\usepackage{centernot}
\usepackage{comment}

\usepackage{hyperref}
\hypersetup{
colorlinks,
urlcolor=azul,
linkcolor=azul,
citecolor=azul,
} 

\definecolor{azul}{RGB}{33,64,154}

\newtheorem{theorem}{Theorem}

\newtheorem{definition}{Definition}
\newtheorem{remark}{Remark}

\numberwithin{equation}{section}

\title[H{\"o}lder regularity for fully nonlinear nonlocal equations]{H{\"o}lder regularity for fully nonlinear nonlocal \\ equations with gradient terms}

\author[J.P. Cabeza]{J. Pablo Cabeza}

\newcommand{\AuthorInfo}{%
  \bigskip\hrule\bigskip
  \noindent
  \footnotesize
  \textsc{\noindent Departamento de Ingeniería Matemática, Universidad de Chile}\\
  \noindent Beauchef 851, Torre Norte, Santiago, Chile\\
  \textit{E-mail:} \href{mailto:jcabeza@dim.uchile.cl}{jcabeza@dim.uchile.cl}%
}

\begin{document}


\begin{abstract}
In this paper we prove H\"older regularity results for viscosity solutions of fully nonlinear nonlocal uniformly elliptic second order differential equations with local gradient terms. This extends the nonlocal counterpart of the work of G. Barles, E. Chasseigne and C. Imbert in JEMS, 2011, to fully nonlinear extremal nonlocal operators. 
\end{abstract}  

\date{October, 2024}

\maketitle


\section{Introduction}

In this chapter, we study interior H\"older regularity of viscosity solutions for nonlocal equations of the form,
\begin{equation} \label{eq} 
\mathcal{F}u(x) + b(x)|Du|^{m} = f(x)\quad\text{in}\quad \Omega, 
\end{equation}
where $u$ is a real valued function, $Du$ stands for the gradient of $u$, and $\mathcal{F}$ is an integral extremal fully nonlinear operator of Pucci or Isaacs type. Also, $b \in C^{0,\tau}$ with $\tau\in(0,1)$, $m\in(0,2]$, and $\Omega \subset \mathbb{R}^n$ is a bounded domain with $n\geq 1$.

We are interested in studying the regularity of viscosity solutions to equations of the form \eqref{eq}, under minor hypotheses about $b$ and $f$.

Over the years, significant attention has been devoted to the study of integro-differential operators and fully nonlinear operators characterized by nonlinear growth in the gradient, largely attributed to their applications. 

In general, two approaches are commonly used to study H\"older estimates for fully nonlinear equations. The first approach, based on the Harnack inequality, can be seen in the work of L. Caffarelli and L. Silvestre \cite{caffarelli2009regularity} and considers equations arising from stochastic control problems with Levy processes. The second approach is a method related to the theory of viscosity solutions introduced by Ishii \& Lions \cite{ishii1990viscosity}, initially designed to get comparison principles, that yields results for second-order nonlinear equations, possibly degenerate, but with continuous coefficients. Furthermore, H\"older estimates were established in \cite{bass2005holder} via probability techniques, and Silvestre in \cite{silvestre2006holder} proved H\"older continuity for a class of differential equations associated to jump processes, via growth lemmas.

As far as local drift terms are concerned, in the work \cite{barles2010holder} by G. Barles, E. Chasseigne, C. Imbert, mixed-type operators are studied, where the local part can be fully nonlinear, but the nonlocal part is linear. In this work, we will focus on exploring the nonlocal part of the operator, i.e., focusing on the nonlocal character of the equation. We will show that the methods used in \cite{bass2005holder} can be extended by combining the techniques of \cite{andrade2024regularity} to produce regularity results for Pucci operators.

We note that H\"older regularity provide a lot of applications, including regularity $C^{1,\alpha}$ obtained by L. Caffareli \& L. Silvestre in \cite{caffarelli2009regularity} and existence of eigenvalues for fully nonlinear integro-differential elliptic equations with a drift term \cite{quaas2020principal}, just to mention a few. See also the parabolic counterparts with gradient growth equal to one in \cite{silvestre2012differentiability}.

Our main result reads as follows
\begin{theorem}\label{theorem1}
    Let $u\in C\left(\overline{\Omega}\right)$ be a viscosity solution of the equation \eqref{eq}, assuming that $s\in(0,1)$, $b\in C^{0,\tau}\left(\overline{\Omega}\right)$ with $\tau\in[0,1)$ and $f\in L^{\infty}\left(\Omega\right)$. If $m\in(0,1)$, then $u$ is $\alpha-$H\"older continuous with exponent
    \begin{equation*}
        \alpha\leq\frac{2s-m+\tau}{1-m}.
    \end{equation*} 
    If $m\in(1,2)$, then $u$ is $\alpha-$H\"older continuous with exponent 
    \begin{equation*}
        \alpha\leq\frac{2s-m+1}{2-m}.
    \end{equation*}
    If $m=2$, then $u$ is $\alpha-$H\"older continuos for any $\alpha\in(0,1)$.
\end{theorem}

\vspace{1em}

\begin{remark}
If $m \in (0,1]$ the conclusion of the theorem holds for any $s\geq 1/2$ with the only assumption of continuity over $b$; whereas it holds for any $s$ if we assume that $b$ is $C^{0,\tau}$ with $\tau\geq 1-2s$. If $m \in (1,2)$ the conclusion of the theorem holds for any $s\geq 1/2$ with merely continuity over $b$; whereas it holds for any $s$ if we assume that $b$ is $C^{0,\tau}$ with $\tau\geq 1-2s$. Note that if $m=2$, the conclusion holds for any $s\in(0,1)$ with the only assumption of the continuity over $b$.
\end{remark}

The rest of the paper follows this structure: the second section presents assumptions about the operators under study and definitions. The third section is dedicated to proving the main theorem in which we make use of the Ishii-Lion's method.

\section{Preliminaries}

In this section, we will provide the basic concepts and assumptions used in this paper. First, we present the definition of the operators we will work in this article, followed by the definition of viscosity solution.

A linear nonlocal uniformly elliptic operator is given by
\begin{align}\label{eq1}
    L_Ku(x) &= \text{P.V. }\int_{\mathbb{R}^n}(u(x+z)-u(x))K(z)\mathrm{d}z \\
    &= \lim_{\epsilon\to0^+}\int_{B_{\epsilon}^c(x)}(u(x+z)-u(x))K(z)\mathrm{d}z \nonumber
\end{align}
where $\text{P.V. }$ denotes the principal value, and the kernel $K:\mathbb{R}^{n}\to\mathbb{R}$ satisfies:
\begin{enumerate}
    \item[(K.1)] is nonnegative, symmetric and 
    \begin{equation*}\label{K1}
        \displaystyle\int_{\mathbb{R}^{n}}\min\{|z|^{2},1\}K(z)\mathrm{d}z < \infty;
    \end{equation*}
    \item[(K.2)] for some $s\in(0,1)$ and ellipticity constants $0<\lambda\leq\Lambda$ we have
    \begin{equation*}\label{K2}
    \dfrac{\lambda}{|z|^{n+2s}}\leq K(z)\leq \dfrac{\Lambda}{|z|^{n+2s}}.
    \end{equation*}
\end{enumerate}

The most classical example of a nonlocal linear uniformly elliptic operator is provided by the fractional Laplacian by considering a kernel of the form $|z|^{-n-2s}$, that is,
\begin{equation*}
    (-\Delta)^{s}u(x) := C_{n,s}\text{ P.V. }\int_{\mathbb{R}^n}\dfrac{u(x)-u(y)}{|x-y|^{n+2s}}\mathrm{d}y \quad ,\forall x\in\mathbb{R}^n,
\end{equation*}
where $C_{n,s}>0$ is a normalization constant. We observe that, via
\begin{equation*}
    \int_{\mathbb{R}^n}\dfrac{u(x+z)-u(x)}{|z|^{n+2s}}\mathrm{d}z = \int_{\mathbb{R}^n}\dfrac{u(x-z)-u(x)}{|z|^{n+2s}}\mathrm{d}z,
\end{equation*}
so the fractional Laplacian can be also written as
\begin{equation*}
    -(-\Delta)^{s}u(x) = \frac{1}{2}C_{n,s}\int_{\mathbb{R}^n}\dfrac{u(x+z) + u(x-z) - 2u(x)}{|z|^{n+2s}}\mathrm{d}z.
\end{equation*}

Let $\mathcal{L}_0$ be the class of all linear uniformly elliptic operators. In other words, we say that $L\in\mathcal{L}_0$ if the operator $L$ has the form \eqref{eq1}, under the same above assumptions about $(\text{K.1})$ and $(\text{K.2})$ over $K$. We denote the nonlocal Pucci operators as
\begin{equation*}
    \mathcal{M}^{+}u(x)=\sup_{L\in\mathcal{L}_0} Lu(x) \quad,\quad\mathcal{M}^{-}u(x)=\inf_{L\in\mathcal{L}_0} Lu(x)
\end{equation*}
i.e., the maximum and minimum operators for the class $\mathcal{L}_0$, respectively. It is noteworthy that, $\mathcal{L}_0$, $\mathcal{M}^+$ and $\mathcal{M}^-$ depend on the parameters $\Lambda$, $\lambda$ and $s$, in particular 
\begin{align*}
    &\mathcal{M}^{+}u(x)=\int_{\mathbb{R}^n}\dfrac{S_{+}(u(x+z) + u(x-z) - 2u(x))}{|z|^{n+2s}}\mathrm{d}z, \\
    &\mathcal{M}^{-}u(x)=\int_{\mathbb{R}^n}\dfrac{S_{-}(u(x+z) + u(x-z) - 2u(x))}{|z|^{n+2s}}\mathrm{d}z,
\end{align*}
where $S_{+}(t)=\Lambda t^{+}-\lambda t^{-}$ and $S_{-}(t)=\lambda t^{+}-\Lambda t^{-}$. Here, $t^{+}$ and $t^{-}$ are defined as $t^{+} = \max\{t, 0\}$ and $t^{-} = \max\{-t, 0\}$, respectively.

Throughout this paper, we will fix the maximal Pucci extremal operator of the form
\begin{equation*}
    \mathcal{M}^{+}u(x)=\sup_{L\in\mathcal{L}_0} Lu(x)
\end{equation*}
where the operator $L$ is of the form, 
\begin{equation*}
    \displaystyle\int_{\mathbb{R}^{n}}(u(x+z)-u(x))K(z)\mathrm{d}z,
\end{equation*}
and the kernel $K:\mathbb{R}^{n}\backslash\{0\}\to\mathbb{R}$ represents the frequency of jumps in every direction $z\in\mathbb{R}^n$ and satisfies some conditions.

In proving H\"older regularity for the main problem, a suitable notion of a viscosity solution has to be considered. We recall the definition of viscosity solution for \eqref{eq} as in \cite[Definition 2.1]{barles2010holder}.

\begin{definition} \label{definition1}
        We say that $u:\mathbb{R}^{n}\to\mathbb{R}$ is a viscosity subsolution to~\eqref{eq} if, for any test function $\varphi \in C^2\left(B_\delta(x_0)\right)$ such that $u - \varphi$ attains its maximum at $x_0\in\Omega$ for some $\delta>0$, we have
\begin{equation*}
    \mathcal{M}^{+}(u,\varphi,x_0,\delta) + b(x_0)|D\varphi(x_0)|^{m}\geq f(x_0) \text{ in } \Omega
\end{equation*}
where the operator $\mathcal{M}^{+}$ is given by
\begin{equation*}
    \mathcal{M}^{+}(u,\varphi,x,\delta) = \sup_{K}\left\{\int_{\mathbb{R}^{n}\backslash B_{\delta}}(u(x+z)-u(z))K(z)\mathrm{d}z \right.\\
    + \left. \int_{ B_{\delta}}(\varphi(x+z)-\varphi(z))K(z)\mathrm{d}z \right\}.
\end{equation*}
    In a similar manner, we say that $u:\mathbb{R}^{n}\to\mathbb{R}$ is a viscosity supersolution to~\eqref{eq} if, for any test function $\varphi \in C^{2}\left(B_{\delta}(x_0)\right)$ such that $u - \varphi$ attains its minimum at $x_0\in\Omega$ for some $\delta>0$, we have
\begin{equation*}
    \mathcal{M}^{+}(u,\varphi,x_0,\delta) + b(x_0)|D\varphi(x_0)|^{m}\leq f(x_0) \text{ in } \Omega
\end{equation*}

A continous function $u:\mathbb{R}^{n}\to\mathbb{R}$ is a viscosity solution to~\eqref{eq} if it is both a viscosity subsolution and supersolution.
\end{definition}

\section{Proof of Theorem \ref{theorem1}}

This section is devoted to the proof of the main theorem of this paper. For the proof, we will use the method introduced by Ishii-Lions \cite{ishii1990viscosity}. Following this method, we fix $x_{0}\in\Omega$ and introduce the function $\displaystyle\Phi:\mathbb{R}^{n}\times\mathbb{R}^{n}\to\mathbb{R}$ given by
 
\begin{equation*}
    \Phi(x,y) := u(x) - u(y) - \phi(x-y) - \Gamma(x),
\end{equation*}
where $\phi$ is a radial function and $\Gamma(x)=L_2|x-x_0|^{2}$ is a localized function around $x_0 \in\Omega$. Since our goal is to attain H\"older regularity, we will choose the following function $\phi(t)=L_{1}|t|^{\alpha}$, $L_{1}>0$ and $\alpha\in(0,1)$.

Our aim is to show that \begin{equation}\label{proof}
    \Phi(x,y)\leq0.
\end{equation}
The first observation to note is that if $\Phi$ is nonpositive, then under suitable control conditions over $L_1$, $L_2$, and $\alpha$, we obtain the desired regularity result
\begin{equation*}
    u(x)-u(y)\leq L_{1}|x-y|^{\alpha}.
\end{equation*}

To prove \eqref{proof}, we will demonstrate that the maximum is only achieved when $x=y=x_0$. To establish this, from the continuity of the function $\Phi$, we observe that it attains a maximum point $(\overline{x},\overline{y})\in\overline{\Omega}\times\overline{\Omega}$. 

Suppose by contradiction that $\Phi(\overline{x},\overline{y})>0$.
As a consequence of this, it follows that,
\begin{enumerate}
    \item $\overline{x}\neq\overline{y}$, since otherwise $\Phi(\overline{x},\overline{y})=-\Gamma(\overline{x})<0$;
    \item $L_{1}|\overline{x}-\overline{y}|^{\alpha}, L_{2}|\overline{x}-x_0|^{2}\leq u(\overline{x}) - u(\overline{y}) \leq 2\|u\|_{\infty};$
    \item $0<u(\overline{x}) - u(\overline{y})$, this implies that $u(\overline{y}) < u(\overline{x})$;
\end{enumerate} 
In particular, if $L_{2}\geq 8\|u\|_{\infty}[d(x_0,\partial\Omega)]^{-2}$, we obtain
\begin{equation*}
    |\overline{x}-x_0|^{2}\leq 2\frac{\|u\|_{\infty}}{L_2}\leq\frac{1}{4}[d(x_0,\partial\Omega)]^2.
\end{equation*}
From this it follows that $|\overline{x}-x_0| \leq d(x_0,\partial\Omega)/2$. Therefore, if $L_1$ is sufficiently large, such that \begin{equation*}
    \left(2\frac{\|u\|_{\infty}}{L_1}\right)^{1/\alpha} < \frac{1}{2}d(x_0,\partial\Omega),
\end{equation*} 
we can see that $\overline{x}, \overline{y}\in \Omega$.

On the other hand, we have that $\Phi_{\overline{y}}(x) = u (x) -\varphi(x,\overline{y})$ reaches a maximum in $B_{\delta}(\overline{x})$ at the point $\overline{x}$, where $\varphi(x,y) = \phi(x-y) + \Gamma(x)$. Then $\varphi\in C^{2}\left(\Omega\right)$ is a test function in the definition of $u$ as a viscosity subsolution of \eqref{eq}, so it follows
\begin{equation*}
    \mathcal{M}^{+}(u,\varphi, \overline{x}, \delta) + b(\overline{x})|D_{x}\varphi(\overline{x}, \overline{y})|^{m}\geq f(\overline{x}).
\end{equation*}
Similarly, $\Phi_{\overline{x}}(y) = u(y)+\varphi(\overline{x},y)$ reaches a minimum in $B_{\delta}({\overline{y}})$ at the point $\overline{y}$; we obtain the viscosity supersolution inequality,
\begin{equation*}
    \mathcal{M}^{+}(u,-\varphi, \overline{y}, \delta) + b(\overline{y})|-D_{y}\varphi(\overline{x}, \overline{y})|^{m}\leq f(\overline{y}).
\end{equation*}
Combining the previous viscosity inequalities, we obtain the following expression that accounts for both local and nonlocal terms,

\begin{multline}\label{desigualdadpucci}
    \mathcal{M}^{+}(u,\varphi, \overline{x}, \delta) - \mathcal{M}^{+}(u,-\varphi, \overline{y}, \delta) \\ + b(\overline{x})|D_{x}\varphi(\overline{x},\overline{y})|^{m} 
    -b(\overline{y})|D_{y}\varphi(\overline{x},\overline{y})|^{m} - f(\overline{x}) + f(\overline{y})\geq0.
\end{multline}

We need to establish that the right-hand side of the above inequality is negative. To this end, we proceed to estimate both local and nonlocal terms. We denote the local terms in~\eqref{desigualdadpucci} by




\begin{align}\label{estimacionlocal}
& T_{l}:= (b(\overline{x})-b(\overline{y}))|D_{x}\varphi(\overline{x}, \overline{y})|^{m} \nonumber\\
& \qquad \qquad \qquad + b(\overline{y})(|D_{x}\varphi(\overline{x}, \overline{y})|^{m}-|D_{y}\varphi(\overline{x}, \overline{y})|^{m}) \\
& \qquad \qquad \qquad \qquad\quad -(f(\overline{x})+f(\overline{y})), \nonumber
\end{align}
where the derivatives for the function $\varphi(x,y)=\phi(x-y)+\Gamma(x)$ are given by \begin{align*}
    D_{x}\varphi(\overline{x}, \overline{y}) &= \alpha L_{1}|\overline{x}-\overline{y}|^{\alpha-2}(\overline{x}-\overline{y}) + 2L_{2}(\overline{x}-x_0), \\
    D_{y}\varphi(\overline{x}, \overline{y}) &= -\alpha L_{1}|\overline{x}-\overline{y}|^{\alpha-2}(\overline{x}-\overline{y}).
\end{align*}

To estimate the local terms in inequality \eqref{desigualdadpucci}, we must distinguish between two cases for the exponent $m$, in addition to making use of the fact that $b$ belongs to the H\"older space $C^{0,\tau}$, so we can find a constant $C_b>0$ such that $|b(\overline{x})-b(\overline{y})|\leq C_{b}|\overline{x}-\overline{y}|^{\tau}$. When $\tau=0$ we only need to use that $b$ is bounded.

First, let us note that if $m\in (0,1]$, then $|p+q|^{m}-|p|^{m} \leq |q|^{m}$. In this case, it follows from \eqref{estimacionlocal}
\begin{align*}
& T_{l}\leq C_{b}|\overline{x}-\overline{y}|^{\tau}\left[(\alpha L_{1}|\overline{x}-\overline{y}|^{\alpha-1})^{m} +(2L_{2}|\overline{x}-x_0|)^{m}\right] \\
& \qquad\quad + (2L_{2}|\overline{x}-x_0|)^{m}\|b\|_{\infty} + 2\|f\|_{\infty}.
\end{align*}

On the other hand, if $m \geq 1$ we know that for $|p|, |q| \geq 0$, 
\begin{equation*}
    (|p|+|q|)^{m} \leqslant 2^{m-1}\max\{|p|^{m},|q|^{m}\}\leqslant 2^{m-1}(|p|^{m} + |q|^{m}).
\end{equation*} 
Furthermore, the function $t \to |t|^{m}$ is convex for $t\in(0,+\infty)$, this implies $|p+q|^{m} \leq |p|^{m} + m|p+q|^{m-2}\langle p+q, q\rangle$. Using the above inequality, it follows that
\begin{equation*}
    |p+q|^{m} - |p|^m \leqslant m2^{m-1}(|p|^{m-1}+|q|^{m-1})|q|.
\end{equation*}
In this case, of \eqref{estimacionlocal}, we have
\begin{align*}
T_{l} & \leq C_{b}2^{m-1}|\overline{x}-\overline{y}|^{\tau} \bigg[ \Big(\alpha L_1|\overline{x}-\overline{y}|^{\alpha-1}\Big)^{m} + \Big(2L_{2}|\overline{x}-x_0|\Big)^m \bigg] \\
&\qquad\quad+ m2^{m} L_{2}\bigg[\big(\alpha L_{1}|\overline{x}-\overline{y}|^{\alpha-1}\big)^{m-1}\\
&\quad\quad\qquad\quad\qquad\qquad+\big(2L_{2}|\overline{x}-x_0|\big)^{m-1}\bigg]|\overline{x}-x_0|\|b\|_{\infty} + 2\|f\|_{\infty}.
\end{align*}

\noindent Given these inequalities, obtained for the local term when $m\in(0,1]$ and $m > 1$, note that $\Phi(\overline{x},\overline{y})>0$ implies that $|\overline{x}-x_0|< 2\|u\|_{\infty}$, therefore, we can choose $L_{1}$ to be sufficiently large in order to control the terms that depend on $\alpha$, $L_{2}$, $\|b\|_{\infty}$, $\|f\|_{\infty}$.


Now we need to estimate the remaining terms of \eqref{desigualdadpucci}, which involve the nonlocal terms. Following \cite{andrade2024regularity}, for given $\epsilon>0$ and $\bar{x} \in \Omega$ fixed, we have \begin{equation*}
    \mathcal{M}^{+}u(\overline{x})\leq L_{K}u(\overline{x})+\epsilon
\end{equation*} for some $L_{K}\in \mathcal{L}_0$. Indeed, this follows from the definition of supremum, since $\mathcal{M}^{+}u(\overline{x})-\epsilon$ is not an upper bound for the set $\{L_{K}u(\overline{x})\}$, since it is less than $\mathcal{M}^{+}u(\overline{x})=\sup L_{K}u(\overline{x})$. Therefore, we can assert that there exists a kernel $K$ such that 
\begin{equation*}
    \mathcal{M}^{+}u(\overline{x})\leq L_{K}u(\overline{x})+\epsilon \qquad\text{ and }\qquad \mathcal{M}^{+}u(\overline{y})\geq L_{K}u(\overline{y}).
\end{equation*} 
In other words, we only need to consider such a linear operator in the rest of the argument,
\begin{equation*}
    L_{K}u(\overline{x}) - L_{K}u(\overline{y}) + \epsilon \geq \mathcal{M}^{+}u(\overline{x}) - \mathcal{M}^{+}u(\overline{y}),
\end{equation*}
where the operators $L_{K}$ are of the form \eqref{eq1}.

Given that these integrals may be singular at the origin, we will partitionate the domain of integration into three parts: $A_{1}:=\mathbb{R}^{n}\backslash B_1$, which relates to the nonlocal term $T^{1}_{nl}$; $A_{2}:= \{z\in B_\delta: (1-\eta)|z||\overline{x}-\overline{y}| \leq |\langle z, \overline{x}-\overline{y}\rangle|\}$, with $\eta\in(0,1)$ small enough, is related to the nonlocal term $T^{2}_{nl}$; and $A_{3} := B_1\backslash  A_2$, which relates to the nonlocal term $T^{3}_{nl}$, so that
\begin{equation*}
    L_{K}u(\overline{x}) - L_{K}u(\overline{y}) = T^{1}_{nl} + T^{2}_{nl} + T^{3}_{nl},
\end{equation*}where the nonlocal terms for $i\in\{1,2,3\}$ are of the form, 
\begin{equation*}
    T^{i}_{nl}(\overline{x},\overline{y}) = \int_{A_i} (u(\overline{x}+z)-u(\overline{x}))K(z)\mathrm{d}z 
    - \int_{A_i} (u(\overline{y}+z)-u(\overline{y}))K(z)\mathrm{d}z,
\end{equation*}
where $A_i$ with $i\in\{1,2,3\}$ describe each of the preceding domains. To avoid overloading the notation, we will denote $a := \overline{x}-\overline{y}$. 


For the initial term $T^{1}_{nl}$, there is no indeterminacy away from the origin. Consequently, we deduce that
\begin{equation}\label{estimaciont1}
    T^{1}_{nl}\leq 4\|u\|_{\infty}\int_{\mathbb{R}^{n}\backslash B_1}K(z)\mathrm{d}z\leq C_1.
\end{equation}
where \begin{equation*}
    C_1 := \Lambda\frac{2\sigma(S_{n-1})}{s}\|u\|_{\infty},
\end{equation*}
with $\Lambda>0$ the ellipticity constant and $\sigma(S_{n-1})$ the surface area of the sphere of dimension $n-1$.

Before proceeding with the estimation to evaluate $T^{2}_{nl}$, let us note that since $(\overline{x},\overline{y})$ is the maximum point for $\Phi$, it follows that for any $d,d'\in\mathbb{R}^{n}$, 
\begin{equation*}
    \Phi(\overline{x}+d,\overline{y}+d')\leq\Phi(\overline{x},\overline{y}).
\end{equation*}
More precisely, we have
\begin{align}\label{key}
    \displaystyle\left(u(\overline{x}+d) - u(\overline{x})\right) -& \left(u(\overline{y}+d') - u(\overline{y})\right) \nonumber\\
    &\qquad \leq \phi(a+(d-d')) - \phi(a) \\ 
    &\qquad\qquad \quad+\Gamma(\overline{x}+d) - \Gamma(\overline{x}). \nonumber
\end{align}
Applying this inequality, first with $d=z, d'=0$ and then with $d=0, d'=z$, we obtain
\begin{align*}
    T^{2}_{nl} &\leq \int_{ A_2}\left(\phi(a+z) - \phi(a)\right)K(z)\mathrm{d}z \\
    &\qquad + \int_{ A_2}\left(\Gamma(\overline{x}+z) - \Gamma(\overline{x})\right)K(z)\mathrm{d}z \\
    &\qquad\quad + \int_{ A_2}\left(\phi(a-z) - \phi(a)\right)K(z)\mathrm{d}z.
\end{align*}

Now, applying Taylor's formula with an integral remainder (for instance Theorem 9.9 in \cite{leoni2023first}) to the function $\phi(a+z)$
    \begin{equation*}
    \phi(a+z) = \phi(a) + D\phi(a)\cdot z + \int_{0}^{1}(1-s)\langle D^{2}\phi(a+sz)\cdot z,z\rangle\mathrm{d}s.
    \end{equation*}
Analogously, we can derive a similar expression for $\phi(a-z)$,
\begin{equation*}
    \phi(a-z) = \phi(a) - D\phi(a)\cdot z + \int_{-1}^{0}(1+s)\langle D^{2}\phi(a+sz)\cdot z, z\rangle\mathrm{d}s,
\end{equation*}
adding these two equations, we observe that the linear terms in $z$ cancel, leaving us with the following result
\begin{equation*}
    \phi(a+z)+\phi(a-z)-2\phi(a) = \int_{-1}^{1}(1-|s|)\langle D^{2}\phi(a+sz)z,z\rangle\mathrm{d}s.
\end{equation*}

Observe that the derivative of the function $\phi(t)=L_{1}|t|^{\alpha}$ is given by
\begin{equation*}
    D\phi(a) = \alpha L_{1}|a|^{\alpha-2}a.
\end{equation*}

For the second derivative, denoting by $\widehat{a}=a/|a|$,
\begin{align*}
    D^{2}\phi(a) &= \alpha L_{1}|a|^{\alpha-2}\left((\alpha-2)\frac{a}{|a|}\otimes\frac{a}{|a|}+I_n\right) \\
    &= \alpha L_{1}|a|^{\alpha-2}\left((\alpha-2)\widehat{a}\otimes\widehat{a}+I_n\right).
\end{align*}
In particular,
\begin{align*}
    D^{2}\phi(a+sz)\cdot z = \alpha L_{1}|a+sz|^{\alpha-2}\left((\alpha-2)\widehat{(a+sz)}\otimes\widehat{(a+sz)}\cdot z + zI_{n}\right).
\end{align*}
As we seek to estimate the Taylor expansion, we observe that

\begin{align}\label{segundaderivada}
    \langle D^{2}\phi(a+sz) \cdot z,z\rangle &= \alpha L_{1}|a+sz|^{\alpha-2}\left((\alpha-2)\widehat{(a+sz)}\otimes\widehat{(a+sz)}\cdot z + zI_n\right)\cdot z \nonumber \\
    &= \alpha L_{1}|a+sz|^{\alpha-2}\left(|z|^2+(\alpha-2)|\langle \widehat{a+s z}, z\rangle|^2 \right).
\end{align}

Note that within the set $ A_2$, which includes $z\in B_{\delta}$, the following bounds are established with $\delta=|a|\rho_0$ and $\rho_0\in(0,1)$
\begin{align}\label{desiguladades}
|a+sz| &\leq|a|+|s||z| \leq |a|+\delta = |a|\left(1+\rho_0\right) \\
|a+sz| &\geq |a|-|s||z|\geq |a|-\delta = |a|\left(1-\rho_0\right) \nonumber \\
|(a+sz)\cdot z| &= |a \cdot z+s|z|^2| \geq |a\cdot z|-\delta|z| \nonumber \\
& \quad \geq (1-\eta)|a||z|-\rho_0|a||z|  \nonumber \\
& \quad = \left(1-\eta-\rho_0\right)|a||z|. \nonumber 
\end{align}
Through the preceding inequalities and the Taylor expansion \eqref{segundaderivada}, we have
\begin{align*}
&\left\langle D^2 \phi(a+s z) \cdot z, z\right\rangle \\
&\quad\leq \alpha L_{1}|a + sz|^{\alpha-4}\left(|a + sz|^{2}|z|^2+(\alpha-2) |(a+sz)\cdot z|^{2} \right) \\
&\quad\leq \alpha L_{1}|a+s z|^{\alpha-4}\left(|a+sz|^2 + (\alpha-2)\left(1-\eta-\rho_0\right)^2|a|^2\right)|z|^2 \\
&\quad\leq \alpha L_{1}|a|^{\alpha-2}(1+\rho_0)^{\alpha-4}\left((1+\rho_0)^2 -(2-\alpha)(1-\eta-\rho_0)^2\right)|z|^2.
\end{align*}
From this estimate, we define \begin{equation*}
    C_{2} := (1+\rho_0)^{\alpha-4}\left((2-\alpha)(1-\eta-\rho_0)^{2}-(1+\rho_0)^{2}\right)
\end{equation*} then,
\begin{equation*}
    \int_{ A_2}\left(\phi(a+z) + \phi(a-z) - 2\phi(a)\right)K(z)\mathrm{d}z \leq -C_{2}\alpha L_{1}|a|^{\alpha-2}\int_{ A_2}|z|^{2}K(z)\mathrm{d}z.
\end{equation*}

Given $ A_2\subset B_1$, we can observe that this set is indeed a cone, so if $z\in A_2$ then $-z\in A_2$. Following the definition of $\Gamma$,
\begin{equation*}
    \int_{ A_2}\left(\Gamma(\overline{x}+z) - \Gamma(\overline{x})\right)K(z)\mathrm{d}z = L_2\int_{ A_2}(|z|^{2}+\langle \overline{x}-x_0,z\rangle)K(z)\mathrm{d}z,
\end{equation*}
by using the symmetry of the kernel, we obtain \begin{equation*}
    L_2\int_{ A_2}(2\langle \overline{x}-x_0,z\rangle)K(z)\mathrm{d}z = 0,
\end{equation*}
that is,
\begin{equation*}
    \int_{ A_2}\left(\Gamma(\overline{x}+z) - \Gamma(\overline{x})\right)K(z)\mathrm{d}z = L_2\int_{ A_2}|z|^{2}K(z)\mathrm{d}z.
\end{equation*}

Given that we are in the \(A_2\), we need to estimate the integral of the kernel. For this, note that
\begin{align*} 
\int_{ A_2}|z|^2 K(z)\mathrm{d}z \leq \sigma(S_{n-1})\Lambda \int_{0}^{\delta} r^{1-2s}\mathrm{d}r = \Lambda \frac{\sigma(S_{n-1})}{2(1-s)}\delta^{2-2s}.
\end{align*}
Hence, by defining this constant
\begin{equation*}
    C'_{2} := \Lambda \frac{\sigma(S_{n-1})}{2-2s}\rho_0^{2-2s},
\end{equation*}
we can see that the second nonlocal term $T^{2}_{nl}$ is bounded by
\begin{equation}\label{estimacion2}
    T^{2}_{nl}\leq -C_{2}C'_{2} \alpha L_{1}|a|^{\alpha - 2s} + L_2C'_{2}|a|^{2-2s}.
\end{equation}

In order to estimate the final nonlocal term $T^{3}_{nl}$ in $B_1\backslash A_2$, that is,
\begin{align*}
    T^{3}_{nl} &\leq \int_{B_1\backslash A_2}\left(\phi(a+z) - \phi(a)\right)K(z)\mathrm{d}z \\
    &\qquad + \int_{B_1\backslash A_2}\left(\Gamma(\overline{x}+z) - \Gamma(\overline{x})\right)K(z)\mathrm{d}z \\
    &\qquad\quad + \int_{B_1\backslash A_2}\left(\phi(a - z) - \phi(a)\right)K(z)\mathrm{d}z 
\end{align*}
we will write $B_1\backslash A_2 = (B_1\backslash B_\delta)\cup (B_{\delta}\backslash A_2)$. In this way, in the annulus $B_1\backslash B_\delta$, we can use a Taylor expansion for $\phi$, while in the domain $B_\delta\backslash A_2$, we utilize the concavity of the function $t^{\alpha}$ for $t\in(0,\infty)$.

Let us consider $d\in\mathbb{R}^n$, then in $B_1\backslash B_\delta$
\begin{align*}
    \phi(a+d) - \phi(a) &\leq L_{1}(|a|+|d|)^\alpha - L_{1}(|a|)^\alpha \\
    &\leq \alpha L_{1}|a|^{\alpha-1}|d|.
\end{align*}
In this way, taking first $d=z$ and then $d=-z$, we can conclude that in $B_1\backslash B_{\delta}$, 
\begin{align*}
\int_{B_1\backslash B_\delta} (\phi( a +z)+\phi( a -z) -& 2\phi( a ))K(z)\mathrm{d}z \\
&\leq 2\alpha L_{1}| a |^{\alpha-1}\int_{B_1\backslash B_\delta} |z|K(z)\mathrm{d}z.
\end{align*}

Given the last integral with $s\in(0,1)$, we must separate the analysis in some cases. If $s\neq 1/2$, then
\begin{equation*}
\int_{B_1\backslash B_\delta} |z|K(z)\mathrm{d}z \leq \Lambda \sigma(S_{n-1})\int_{\delta}^{1}r^{-2s}\mathrm{d}r = \Lambda \frac{\sigma(S_{n-1})}{1-2s}\left(1-\delta^{1-2s}\right),
\end{equation*}
and we must again distinguish into two cases, as if $s<1/2$
\begin{equation*}
\int_{B_1\backslash B_\delta} |z|K(z)\mathrm{d}z \leq \Lambda \frac{\sigma(S_{n-1})}{1-2s}=\overline{C}_3,
\end{equation*}
whereas if $s>1/2$ with $\delta=|a|\rho_0$
\begin{equation*}
\int_{B_1\backslash B_\delta} |z|K(z)\mathrm{d}z \leq \Lambda \frac{\sigma(S_{n-1})}{2s-1}\left(|a|\rho_0\right)^{1-2s}=|a|^{1-2s}\widetilde{C}_3.
\end{equation*}
If we define the following constant $C_{3} = \max\{\overline{C}_3, \widetilde{C}_3\}$ we obtain the following upper bound for $\phi$ in $B_1\backslash B_\delta$, when $s<1/2$,
\begin{align*}
\int_{B_1\backslash B_\delta} (\phi( a +z)+\phi( a -z) -& 2\phi( a )) K(z)\mathrm{d}z \leq 2\alpha L_{1}C_{3}|a|^{\alpha-1},
\end{align*}
and if $s>1/2$,
\begin{align*}
\int_{B_1\backslash B_\delta} (\phi( a +z)+\phi( a -z) -& 2\phi( a )) K(z)\mathrm{d}z \leq 2\alpha L_{1}C_{3}| a |^{\alpha-2s}.
\end{align*}

In the case $s=1/2$ we get
\begin{equation*}
\int_{B_1\backslash B_\delta} |z|K(z)\mathrm{d}z \leq -\Lambda \sigma(S_{n-1})\ln(|a|\rho_0),
\end{equation*}
and if we define the constant \begin{equation*}
    C'_{3} := \Lambda\sigma(S_{n-1})
\end{equation*}
we obtain the following upper bound,
\begin{align*}
\int_{B_1\backslash B_\delta} \phi( a +z)+\phi( a - z) -& 2\phi( a )) K(z)\mathrm{d}z \\
&\leq -2\alpha L_{1}C'_{3}| a |^{\alpha-1}\ln(|a|\rho_0).
\end{align*}

In $B_\delta\backslash A_2$ we use a Taylor expansion, since $\alpha\in(0,1)$ we see that
\begin{equation*}
    \phi(a+z)+\phi(a-z)-2 \phi(a) = \int_{-1}^{1}(1-|s|)\langle D^{2}\phi(a+sz)z,z\rangle\mathrm{d}s.
\end{equation*} 
In particular,
\begin{align*}
\langle D^{2}\phi(a+sz) \cdot z,z\rangle &= \alpha L_{1}|a+sz|^{\alpha-2}\left(|z|^2-(2-\alpha)|\langle \widehat{a+s z}, z\rangle|^2 \right) \\
&\leq \alpha L_{1}|a+sz|^{\alpha-2}|z|^2.
\end{align*}
From \eqref{desiguladades}, we can conclude that
\begin{align*}
\langle D^{2}\phi(a+sz) \cdot z,z\rangle \leq \alpha L_{1}|a|^{\alpha-2}(1-\rho_0)^{\alpha-2}|z|^2,
\end{align*}
since the kernel $K$ is bounded by constant $\Lambda>0$, the estimate of the integral over $B_{\delta}\backslash A_2$ is given by
\begin{align*}
    \int_{B_\delta\backslash A_2}\frac{1}{|z|^{n+2s-2}}\mathrm{d}z & \leq \Lambda\sigma(S_{n-1})\int_{0}^{\delta}\frac{1}{r^{2s-1}}\mathrm{d}r = \Lambda\frac{\sigma(S_{n-1})}{2(1-s)}\delta^{2-2s}.
\end{align*} 
Now, if we define \begin{equation*}
    C_4 := \Lambda\frac{\sigma(S_{n-1})}{2(1-s)}\rho_0^{2-2s},
\end{equation*}
we obtain the following bound for $\phi$,
\begin{align*}
\int_{B_{\delta}\backslash  A_2} \left(\phi( a +z)+\phi( a -z) -2\phi( a )\right)K(z)\mathrm{d}z  \leq \alpha(1-\rho_0)^{\alpha-2}L_{1}C_4| a |^{\alpha-2s}.
\end{align*}

Finally, for the integral involving $\Gamma$, we can observe that
\begin{align*}          
\int_{B_1\backslash A_2}\left(\Gamma(\overline{x}+z) - \Gamma(\overline{x})\right)K(z)\mathrm{d}z 
&\leq \int_{B_1}\left(\Gamma(\overline{x}+z) - \Gamma(\overline{x})\right)K(z)\mathrm{d}z \\
&\leq L_{2}\int_{B_1}|z|^{2}K(z)\mathrm{d}z.
\end{align*}
In a similar manner to the previous kernel estimation, 
\begin{align*}
    \int_{B_1}\frac{1}{|z|^{n+2s-2}}\mathrm{d}z \leq \Lambda\frac{\sigma(S_{n-1})}{2(1-s)}=: C'_4,
\end{align*} 
we can conclude that the integral of $\Gamma$ over $B_1\backslash  A_2$ is bounded by
\begin{equation*}
\int_{B_1\backslash A_2}\left(\Gamma(\overline{x}+z) - \Gamma(\overline{x})\right)K(z)\mathrm{d}z \leq L_{2}C'_{4}.
\end{equation*}

In conclusion, when combining all terms related to the nonlocal term $T^{3}_{nl}$, we must take into account the possible values of $s\in(0,1)$. If $s<1/2$,
\begin{align}\label{estimacion<1/2}
T^{3}_{nl} &\leq 2\alpha L_1C_{3}| a |^{\alpha-1} + \alpha(1-\rho_0)^{\alpha-2} L_{1}C_4| a |^{\alpha-2s} + L_{2}C'_{4}, 
\end{align}
if $s=1/2$, 
\begin{align*}
T^{3}_{nl} &\leq -2\alpha L_1\ln(|a|\rho_0)C'_{3}| a |^{\alpha-1} + \alpha(1-\rho_0)^{\alpha-2} L_{1}C_4| a |^{\alpha-1} + L_{2}C'_{4}.
\end{align*}
and if $s>1/2$
\begin{align}\label{estimacion>1/2}
T^{3}_{nl} &\leq 2\alpha L_1 C_{3}| a |^{\alpha-2s} + \alpha(1-\rho_0)^{\alpha-2} L_{1}C_4| a |^{\alpha-2s} + L_{2}C'_{4}, 
\end{align}

In the case $s=1/2$, we observe that $\displaystyle\lim_{r \to 0}\ln(|r|)r^{\alpha}\to 0$. Therefore, we also obtain powers of $|a|$ of the form $\alpha-2s$. Thus, \eqref{estimacion>1/2} remains valid in this case as well.

\subsection{Estimates for local and nonlocal terms}




We can now proceed to a joint analysis of the local and nonlocal contributions appearing in equation \eqref{desigualdadpucci}. To this end, we begin by considering $s\geq1/2$ and examine the result of combining the local terms \eqref{estimacionlocal} and the nonlocal terms \eqref{estimaciont1}, \eqref{estimacion2}, \eqref{estimacion>1/2}, for both cases $m\in(0,1]$ and $m\in(1,2]$.

For \(m\in(0,1]\), the following estimate is obtained from inequality~\eqref{desigualdadpucci},
\begin{equation}\label{M11}
\begin{split}
0 &\leq \nu_1\big(L_2, \|b\|_\infty, \|u\|_\infty, \|f\|_\infty\big) \\
&\qquad - \alpha L_1 |a|^{\alpha - 2s} \Big( 
C_2 C_2' - (1-\rho_0)^{\alpha-2} C_4 
- C_b (\alpha L_1 |a|^\alpha)^{m-1} |a|^{\tau + 2s - m} - 2 C_3 |a|^{2s-1} 
\Big),
\end{split}
\end{equation}
where
\[
    \nu_1(L_2,\|b\|_\infty,\|u\|_\infty,\|f\|_\infty) :=
    \begin{aligned}[t]
        & 2 \|f\|_\infty 
        + C_b |a|^\tau (2 L_2 |\overline{x}-x_0|)^m 
        + \|b\|_\infty (2 L_2 |\overline{x}-x_0|)^m \\
        &\quad + C_1 \|u\|_\infty 
        + L_2 C_2' |a|^{2-2s} 
        + L_2 C_4' 
        + \epsilon.
    \end{aligned}
    \]
The term $\nu_1$ gathers all contributions that are negligible for small $|a|$ or remain bounded independently of $L_1$. In order to conclude, it is sufficient to ensure that





$$(\alpha L_{1}|a|^{\alpha})^{m-1}|a|^{\tau+(2s-m)}=o(1),$$ which holds if $\tau+2s-m+\alpha(m-1)\geq0$ and $m\in(0,1)$, that is, 
\begin{align*}
    \alpha\leq\frac{\tau-m+2s}{1-m}\quad\text{ if }m\in(0,1)\quad\text{ and }\quad\tau+2s-1\geq0\quad\text{ if } m=1.
\end{align*}
We choose $L_1$ sufficiently large such that $L_{1}>2^{1+\alpha}d(x_0,\Omega)^{-\alpha}\|u\|_{\infty}$ and
\begin{equation*}
    L_{1} \geq \frac{\nu_1\big(L_2, \|b\|_\infty, \|u\|_\infty, \|f\|_\infty\big)}{(C_{2}C^{'}_{2} - (1-\rho_{0})^{\alpha-2}C_{4})\alpha}.
\end{equation*}
Then, based on the condition $L_{1}|a|^{\alpha} \leq 2\|u\|_{\infty}$, for $\alpha\in (0, 1)$, we deduce that $|a|^{\alpha - 2s} > 1$ since these terms tend to $0$ as $\overline{x}$ approaches $\overline{y}$. It follows from
\eqref{M11} that
\begin{equation*}
    0\leq \nu_{1}(L_2, \|b\|_{\infty}, f)(1-|a|^{\alpha-2s})<0,
\end{equation*}
which is a contradiction. Therefore, \(\Phi(x, y) \leq 0\).

\noindent Furthermore, in the case when $m=1$, if $s\geq 1/2$ then the conclusion holds for any $\tau\geq0$, that is, $b$ it can be only a continuos function. Now, it follows for $m\in(0,1)$ that
\begin{equation*}
    \alpha\leq\frac{\tau+2s-1}{1-m}.
\end{equation*} 

When considering the case of $m\in(1,2]$, the following estimate is obtained from inequality~\eqref{desigualdadpucci},

\begin{equation}\label{M2}
\begin{split}
0 \le{} &\; \nu_2\big(L_2,\|b\|_\infty,\|u\|_\infty,\|f\|_\infty\big) \\
&\quad - \alpha L_1 |a|^{\alpha-2s}
\Bigg(
C_2 C_2' 
- (1-\rho_0)^{\alpha-2} C_4 
- C_b(\alpha L_1|a|^\alpha)^{m-1}|a|^{\tau+2s-m} \\
&\qquad\qquad\qquad\qquad
- 2 C_3 |a|^{2s-1}
- m2^{m}\|b\|_\infty L_2 (\alpha L_1|a|^{\alpha})^{m-2}|a|^{2s+1-m}|\overline{x}-x_0|
\Bigg),
\end{split}
\end{equation}
where
\[
    \nu_2(L_2,\|b\|_\infty,\|u\|_\infty,\|f\|_\infty) :=
    \begin{aligned}[t]
        & 2\|f\|_\infty 
        + C_b 2^{m-1}|a|^\tau \big[ (\alpha L_1 |a|^{\alpha-1})^m + (2 L_2 |\overline{x}-x_0|)^m \big] \\ 
        & \qquad + m 2^m \|b\|_\infty L_2 (2 L_2 |\overline{x}-x_0|)^{m-1}|\overline{x}-x_0| + C_1 \|u\|_\infty \\
        & \qquad\qquad+ L_2 C_2' |a|^{2-2s} 
        + L_2 C_4' 
        + \epsilon.
    \end{aligned}
    \]
Similarly to the previous case, we need to ensure that 
\begin{equation*}
    |a|^{\tau}(L_{1}|a|^{\alpha})^{(m-1)}|a|^{2s-m}=o(1)\quad\text{ and }\quad |a|(L_{1}|a|^{\alpha})^{m-2}|a|^{2s-m}=o(1),
\end{equation*}
which holds provided that $\tau+\alpha(m-1)+2s-m\geq0$ and $1+\alpha(m-2)+2s-m\geq0$, from where we obtain 
\begin{equation*}
    \alpha\leq\frac{1-m+2s}{2-m}\quad\text{ if }m\in(1,2)\quad\text{ and } \quad\tau +\alpha + 2(s-1)\geq0\quad\text{ if }m=2.
\end{equation*}
We choose $ L_{1}>0 $ large enough to control the terms of \(\nu_2\),
\begin{equation*}
    L_{1}\geq \frac{\nu_{2}(L_{2},\|b\|_{\infty}, \|u\|_{\infty}, \|f\|_{\infty})}{(C_{2}C'_{2}-(1-\rho_0)^{\alpha-2} C_{4})\alpha }.
\end{equation*}%
Since $\alpha\in(0,2s)$, we have $|a|^{\alpha-2s}>1$. It follows that
\begin{equation*}
    0\leq\nu_{2}(L_{2},\|b\|_{\infty}, \|u\|_{\infty}, \|f\|_{\infty})(1-|a|^{\alpha-2s})<0.
\end{equation*}%
This yields the desired contradiction.

Furthermore, in the case when $m=2$, from \eqref{M2}, we impose that $\tau-\alpha+2s\geq0$ and from the restriction $\tau +\alpha + 2(s-1)\geq0$, we obtain that $\tau\geq 1-2s$, however we know that $1-2s\leq0$, so the conclusion holds for any $\tau\geq0$, that is, $b$ it can be only a continuous function.

Now we see what happens when $s<1/2$, using estimates for the local terms \eqref{estimacionlocal} and the nonlocal terms \eqref{estimaciont1}, \eqref{estimacion2}, \eqref{estimacion<1/2} for both cases $m\in(0,1]$ and $m\in(1,2]$. 

In the case $m\in(0,1]$ and $s<1/2$, we obtain the following estimate from inequality~\eqref{desigualdadpucci},
\begin{equation}\label{M3}
\begin{split}
0 \leq{} & 
\nu_3(L_2, \|b\|_\infty, \|u\|_\infty, \|f\|_\infty) \\
&\quad 
- \alpha L_{1}|a|^{\alpha-2s} \Big(C_{2}C^{'}_{2} - (1-\rho_{0})^{\alpha-2}C_{4} - C_{b}|a|^{\tau}(\alpha L_{1}|a|^{\alpha})^{m-1}|a|^{2s-m} - 2C_{3}|a|^{2s-1}\Big),
\end{split}
\end{equation}
where, analogously, we have grouped all terms that are bounded or depend on \(L_2\) into \(\nu_3\).
In order to conclude, it is sufficient to establish that
\begin{equation*}
    |a|^{\tau}(L_{1}|a|^{\alpha})^{(m-1)}|a|^{2s-m}=o(1),
\end{equation*} that is,
\begin{equation*}
    \alpha\leq\frac{\tau-m+2s}{1-m}\quad\text{ if }m\in(0,1)\quad\text{ and } s<\frac{1}{2},
\end{equation*}
while if $m=1$, we obtain $\tau+2s-1\geq0$. In the other words, if $s< 1/2$ and $m\in(0,1]$ then we must have $b$ is $C^{0,\tau}$ with $\tau \geq 1-2s$.

In the case when $m\in(1,2]$, and $s<1/2$, we obtain the estimate from inequality~\eqref{desigualdadpucci}
\begin{equation}\label{M4}
\begin{split}
0 &\le{}  
\nu_4(L_2, \|b\|_\infty, \|u\|_\infty, \|f\|_\infty, |a|) \\
&\qquad - \alpha L_1 |a|^{\alpha-2s} 
\Bigg(
C_2 C'_2 
- 2 C_3 |a|^{2s-1} 
- (1-\rho_0)^{\alpha-2} C_4 \\
&\qquad\qquad\qquad\qquad\qquad- \alpha^{m-1} 2^{m-1} C_b |a|^\tau (L_1 |a|^\alpha)^{m-1} |a|^{2s-m} \\
&\qquad \qquad \qquad \qquad \qquad \qquad- m \alpha^{m-1} 2^m L_2 |a| (L_1 |a|^\alpha)^{m-2} |a|^{2s-m} |\overline{x}-x_0| \|b\|_\infty
\Bigg).
\end{split}
\end{equation}
As before, all terms that are bounded or depend on \(L_2\) have been grouped into \(\nu_4\).
In a completely analogous manner to the previous, we need that
\begin{equation*}
    |a|^{\tau}(L_{1}|a|^{\alpha})^{(m-1)}|a|^{2s-m}=o(1)\quad\text{ and }\quad |a|(L_{1}|a|^{\alpha})^{m-2}|a|^{2s-m}=o(1),
\end{equation*}
that is,
\begin{equation*}
    \alpha\leq \frac{1-m+2s}{2-m}\quad\text{ if }m\in(1,2)\quad\text{ and }\tau+\alpha+2s-2\geq0\quad\text{ if }m=2. 
\end{equation*}

Moreover, in the case when $m=2$, from \eqref{M4}, we impose the condition $\tau-\alpha+2s\geq0$ and from the restriction $\tau+\alpha+2(s-1)\geq0$, we obtain that $\tau\geq1-2s$, however we know that $1-2s>0$. Therefore, the conclusion holds for $b$ in $C^{0,\tau}$ with $\tau\geq1-2s$. 

\vspace{1em}

\begin{remark}
    We mention that the case of a minimal Pucci operator of the form $\mathcal{M}^{-}u(x) = \inf_{L\in\mathcal{L}_0} Lu(x)$ will be analogous, as well as for Isaacs operators in the form
\begin{equation*}
    \mathcal{F}u(x) = \inf_{a}\sup_{b} L_{a,b}u(x)\quad,\quad\mathcal{F}u(x) = \sup_{a}\inf_{b} L_{a,b}u(x)
\end{equation*}
where $L_{a,b}=L$ for some $K=K_{a,b}$ for a family of indexes $a\in A$, $b\in B$, this type of operator is treated as in article \cite{andrade2024regularity}.
\end{remark}

\bigskip

{\small \noindent{\bf Acknowledgments.} The author acknowledges the support of \textit{ANID-Subdirecci\'on de Capital Humano, Doctorado Nacional} 2022-21222157. The author is indebt with Disson dos Prazeres for his hospitality in Sergipe, Brazil, for his help and numerous suggestions during the preparation of this manuscript. Additionally, thank to Gabrielle Nornberg for helpful and extensive discussions on this topic.}

\medskip



\vfill
\AuthorInfo

\end{document}